\documentclass[onefignum,onetabnum]{siamart190516}

\usepackage{lipsum}
\usepackage{amsfonts}
\usepackage{graphicx}
\usepackage{epstopdf}
\usepackage{algorithmic}
\ifpdf
  \DeclareGraphicsExtensions{.eps,.pdf,.png,.jpg}
\else
  \DeclareGraphicsExtensions{.eps}
\fi

\newsiamremark{remark}{Remark}
\newsiamremark{hypothesis}{Hypothesis}
\crefname{hypothesis}{Hypothesis}{Hypotheses}
\newsiamthm{claim}{Claim}

\headers{Maximizing bounded hessian functions}

\title{On global maximization of bounded hessian functions over strongly convex domains}

\author{Costandin Marius\thanks{General Digits 
  (\email{costandinmarius@gmail.com}, \url{http://www.generaldigits.com}).}
}

\usepackage{amsopn}

\usepackage{algorithmic}

\begin{document}

\maketitle

\begin{abstract}
  In this paper we present two frameworks in which global maximization of a bounded hessian function over a strongly convex set can be reduced to convex optimization. The first presented framework is a continuation of one of our previous papers \cite{spheres}. We improve the results and give an explicit algorithm for the computation in polynomial time of the farthest point  in a finite intersection of n-disks to $C \in \mathbb{R}^{n \times 1}$ under the requirement that $C$ does not belong to the convex hull of the centers of the n-disks. Finally, in order to overcome this limitation we present a second framework which characterizes the furthest in the finite intersection of n-disks with $C$ in the convex hull. Unfortunately this second framework requires the ability to decide if a polytope that we define is included in the intersection, which is hard in general. However, as a particular application of our second framework we are able solve in P time some instances of the subset sum problem with real entries: given a set of real numbers decide if there is a subset which adds up to zero.  
\end{abstract}



\section{Introduction}
\subsection{Notations}
We will use throughout the paper the symbol $d(\cdot, \times)$ where $\cdot$ can be a point and $\times$ can be a point or a convex set of points, to designate the Euclidean distance between $\cdot$ and $\times$.  For a vector $u\in\mathbb{R}^n$, $u =\left(u_1, ...,u_n\right)^T$  and $r>0$, we denote by $\mathcal{B}(u,r)$ the open ball centered at $u$ and of radius $r$. We also denote by $\|u\|$, $\|u\|^2 = u^T u$,  the Euclidean norm of the vector $u$. 
With $\frac{\partial}{\partial x}$ we denote the following operator:
\begin{align}
\frac{\partial }{\partial x} = \begin{bmatrix} \frac{\partial }{\partial x_1} &\hdots &\frac{\partial }{\partial x_n}\end{bmatrix}
\end{align}

\subsection{Problem definition}
We are interested in solving the following kind of problems:
\begin{align}
\max_{h(x) \leq 0} f(x)
\end{align} where $h$ is continuous and strongly convex with $\{x | h(x) < 0\} \neq \emptyset$ and $f$ continuous. In the following section, we will add or drop additional requirements about these functions. 

\section{Main results}

Let $f, h : \mathbb{R}^{n \times 1} \to \mathbb{R}$ with $h$ continuous strongly convex, $f$ continuous  such that $k_h \cdot h(x) - k_f \cdot f(x)$ is convex for some $k_h,k_f > 0$. We assume throughout the paper that $\{x | h(x) < 0\} \neq \emptyset$. Next we define the function:
\begin{align}
g(x) = \max \hspace{0.1cm} \{ R - k_f \cdot f(x), 0\} + k_h\cdot h(x)
\end{align} We argue that $g(x)$ is continuous and convex. Indeed, for convexity let us evaluate 
\begin{align}
&g(\lambda \cdot x + (1 - \lambda) \cdot y) = \nonumber \\
& 
\begin{cases}  
R - k_f  f(\lambda \cdot x + (1 - \lambda)y)) +  k_h  h(\lambda  x + (1- \lambda)  y) \hspace{0.2cm} &R - k_f f(\lambda  x + (1 - \lambda) y)) > 0 \\
k_h \cdot h(\lambda \cdot x + (1- \lambda) \cdot y)  &R - k_f f(\lambda x + (1 - \lambda) y)) \leq 0
\end{cases}
\end{align} In the first case, since $k_h \cdot  h - k_f \cdot f$ is convex
\begin{align}
&g(\lambda x + (1 - \lambda) y) \leq \lambda \left( R + k_h h(x) - k_f f(x) \right) + (1 - \lambda) (R + k_h h(y) + k_f f(y)) \leq \nonumber \\
& \leq \lambda (\max \{ R - k_f f(x), 0\} + k_h h(x)) + (1 - \lambda) ( \max \{ R - k_f f(y), 0 \}  + k_h h(y))  \nonumber \\
& \leq \lambda \cdot g(x) + (1 - \lambda) \cdot g(y) 
\end{align} In the second case, since $h$ is convex
\begin{align}
&g(\lambda x + (1 - \lambda) y) \leq \lambda k_h \cdot h (x) + (1 - \lambda) k_h \cdot h(y) \nonumber \\
& \leq \lambda (\max \{ R - k_f f(x), 0\} + k_h h(x)) + (1 - \lambda) (\max \{R - k_f f(y), 0\} + k_h h(y))  \nonumber \\
& \leq \lambda g(x) + (1 - \lambda ) g(y) 
\end{align} 
For continuity we just say that $g$ is the point maximum of continuous functions, hence it is itself continuous. 

The following points are now easy to compute: 
\begin{align}
x^{\star}\in \mathcal{X}^{\star} &= \mathop{\text{argmin}}_{x \in \mathbb{R}^{n \times 1}} \hspace{0.1cm} g(x) \nonumber \\
x_1^{\star} \in \mathcal{X}_1^{\star} &= \mathop{\text{argmin}} \hspace{0.2cm} k_h \cdot h(x) - k_f \cdot f(x) \hspace{0.2cm} \text{s.t} \hspace{0.2cm} h(x) \leq 1
\end{align} 

The following theorem is a central result of this paper:
\begin{theorem} \label{T2.1} If $\{x|h(x) <0 \} \cap \mathcal{X}^{\star}_1 = \emptyset $ then 
\begin{align}
\{x | h(x) < 0\} \hspace{0.2cm} \setminus \{ x | R - k_f \cdot f(x) > 0 \} \neq \emptyset \iff h(x^{\star}) < 0, \hspace{0.2cm} k_f \cdot f(x^{\star}) \geq R
\end{align}
\end{theorem}
\begin{proof} 
The reverse implication being obvious, we focus on the forward implication. We therefore know that $\exists x_0 $ with $k_f \cdot f(x_0) \geq R$ and $h(x_0) < 0$.  Then, it follows that $h(x^{\star}) < 0$ (otherwise $g(x^{\star}) \geq 0 > g(x_0)$ which is impossible since $x^{\star}$ is the minimizer of $g$). The following disjoint reunion is true:
\begin{align}
& \{x | h(x) < 0\} = \nonumber \\
& \left( \{x | h(x) < 0\} \cap \{x | R - k_f \cdot f > 0\}\right) \cup  \left( \{x | h(x) < 0\} \setminus \{x | R - k_f \cdot f > 0\}\right)
\end{align} Therefore $x^{\star}$ is either in the first set, either in the second. We show that the requirements in the theorem statement are sufficient to rule out the first option. Indeed, if $x^{\star} \in \{x | h(x) < 0\} \cap \{x | R - k_f \cdot f > 0\}$ then 
\begin{align}
g(x^{\star}) = R - k_f \cdot f(x^{\star}) + k_h \cdot h(x^{\star}) > g(x_1^{\star}) 
\end{align} Because $x^{\star} \in \{x | h(x) < 0\}$, $x_1^{\star} \in \{ x| h(x) \geq 0\}$, $g(x_1^{\star}) < g(x^{\star})$ and $g$ is continuous and convex, follows that    exists a ball $\mathcal{B}(x^{\star}, \epsilon) \subseteq \{x | h(x) < 0\}$ in the open set $\{x | h(x) < 0\}$ and $y \in \mathcal{B}(x^{\star}, \epsilon)$ such that $g(y) < g(x^{\star})$. Take $y$ on the intersection of the line connecting $x$ and $x^{\star}$ and the sphere $\partial \mathcal{B} \left( x^{\star}, \frac{\epsilon}{2} \right)$.
The existence of $y$ thought should be impossible since $x^{\star}$ is the unique global minimum, hence this is a contradiction.
\end{proof} 

\begin{remark} The required test in Theorem \ref{T2.1} $\{x | h(x) < 0\} \cap \mathcal{X}^{\star}_1 =^{?} \emptyset$ is a convex optimization problem, since it is a feasibility problem of the intersection of two convex sets. 
In the cases where $k_h \cdot h - k_f \cdot f$ is strictly convex, the set $\mathcal{X}^{\star}_1 = \{x^{\star}_1\}$, i.e contains a single point, since the minimizer of strictly convex functions is unique.    
\end{remark}


\begin{remark} It is obvious that Theorem \ref{T2.1} can be used to solve:
\begin{align} \label{E2.10}
\text{argmax} \hspace{0.2cm} f(x) \hspace{0.2cm} \text{s.t} \hspace{0.2cm} h(x) \leq 0
\end{align} for a bounded hessian function $f$ if $h(x_1^{\star}) \geq 0$ for all $x_1^{\star} \in \mathcal{X}^{\star}_1$. Indeed, assume that a point $x_1 \in \{x | h(x) < 0\}$ and an upper bound $\bar{F}$ on the $f(x)$ for $x \in \{x | h(x) < 0\}$ are known. Therefore 
\begin{align}
 \underline{R}  = k_f \cdot f(x_1)  \leq \max_{h(x) \leq 0} k_f \cdot f(x) \leq k_f \cdot \bar{F} = \bar{R} 
\end{align} Take $R = \frac{\underline{R} + \bar{R}}{2}$ and ask $\exists^{?} x \in \{x | h(x) < 0\} \setminus \{x | R - k_f f(x) > 0\}$. This can be answered using Theorem \ref{T2.1} by simply testing if $x^{\star} \in \{x | h(x) < 0\} \setminus \{x | R - k_f f(x) > 0\}$. If the answer is $YES$ then $f(x^{\star}) \geq \frac{R}{k_f}$ and $h(x^{\star}) < 0$ hence let $\underline{R} = k_f f(x^{\star})$, otherwise $\bar{R} = k_f \cdot f(x^{\star})$. This is actually a bisection after the parameter $R$ and will result in the solution of the problem (\ref{E2.10})
\end{remark}

\section{Applications}
In this subsection, we give a few applications to the above theoretical result.
\subsection{Type 1}
Let $f(x), h(x)$ be two continuous and double differentiable  such that $h(x)$ is strongly convex and $k_h \cdot h(x) - k_f \cdot f(x)$ is convex. Then for $C \in \mathbb{R}^{n \times 1}$ of large enough magnitude one has $ k_h \cdot \frac{\partial h}{\partial x} - k_f \cdot \frac{\partial f}{\partial x} + C \neq 0$ for all $x \in \{x | h(x) \leq 0\}$ hence we can apply Theorem \ref{T2.1} to solve the following problem:
\begin{align}
\max \hspace{0.2cm} f(x) + C^T \cdot x \hspace{0.2cm} \text{s.t} \hspace{0.2cm} h(x) \leq 0
\end{align}
\subsection{Type 2}
Next, let us consider the closed n-disks $\bar{\mathcal{B}}(C_k, r_k) = \{x | \|x - C_k\|^2 \leq r_k^2\}$ for some $C_k \in \mathbb{R}^{ n \times 1}$ with $k \in \{1, \hdots, m\}$, $m \in \mathbb{N}$. Let $C \in \mathbb{R}^{n \times 1}$ and we want to solve the following problem:
\begin{align}\label{E3.1}
\max \hspace{0.2cm} \|x - C\|^2 \hspace{0.2cm} \text{s.t} \hspace{0.2cm} x \in \bigcap_{k=1}^m \bar{\mathcal{B}}(C_k, r_k) 
 \end{align}  Let therefore
\begin{align} \label{E3.2}
h(x) &= \max_{k \in \{1, \hdots, m\}} \|x - C_k\|^2 - r_k^2 \nonumber \\
f(x) &= \|x - C\|^2
\end{align} where $h$ is strongly convex, being the maximum of strongly convex functions. We obtain
\begin{align}
g(x) &= \max\{R^2 - \|x - C\|^2, 0\} + \max_{k} \|x - C_k\|^2 - r_k^2 
\end{align} It can be shown that $g$ is continuous and convex. Another way of writing $g$ is
\begin{align}
g &= \max\{ R^2 - \|x - C\|^2 + \max_k \|x - C_k\|^2 - r_k^2, \max_k \|x - C_k\|^2 - r_k\} \nonumber \\
& = \max \{ \max_k \|x - C_k\|^2 - \|x - C\|^2 + R^2 - r_k^2, \max_k \| x - C_k\|^2 - r_k^2\} \nonumber \\
& = \max \{ \max_k 2 \cdot (C - C_k)^T \cdot x + R^2 - r_k^2, \max_k \|x - C_k\|^2 - r_k^2\}
\end{align} that is, the maximum of two functions, each being the point maximum of $m$ convex functions. 

We show that we can solve (\ref{E3.1}) in P time (with the presented theory) if 
\begin{align}\label{E3.6}
C \not\in \text{conv} \{C_1, \hdots, C_m\}
\end{align} Indeed, if this is the case, then $\exists d, H \in \mathbb{R}^{n \times 1}$ such that the hyperplane $\{x | d^T \cdot (x - H) \leq 0\}$ separates $C$ and $\text{conv}\{C_1, \hdots, C_m\}$. This means that 
\begin{align}\label{E3.7}
d^T \cdot (C - H) < 0 \hspace{0.5cm} d^T\cdot (C_k - H) > 0
\end{align} Since $g$ is convex, in order to apply the Theorem \ref{T2.1} it is sufficient to show that $h-f$ does not have a minimum. Indeed, assume that $\exists x_1$ such that $(h-f)(x_1) \leq (h-f)(x)$ for all $x \in \mathbb{R}^{n \times 1}$. Let $x_2 = x_1 + d$ then
\begin{align}
(h - f)(x_2) &= \max_{k} 2\cdot(C - C_k)^T \cdot (x_1 + d) - r_k^2 + \|C_k\|^2 - \|C\|^2\nonumber \\
& \leq \max_k 2 \cdot (C - C_k)^T \cdot x_1 - r_k^2 + \|C_k\|^2 - \|C\|^2 + \max_{k} 2 \cdot (C - C_k)^T \cdot d \nonumber \\
& \leq (h - f)(x_1) + 2 \cdot \max_{k} (C - C_k)^T \cdot d
\end{align} However, using (\ref{E3.7}) it is easy to see that $\max_k (C - C_k)^T \cdot d < 0$ hence $(h - f)(x_2) < (h-f)(x_1)$. This is a contradiction with $x_1$ being the minimizer. 
\begin{remark}
Please note that testing (\ref{E3.6}) boils down to solving a linear feasibility problem. Indeed, we want to know if exists $d \in \mathbb{R}^{n \times 1}$ and $d_H \in \mathbb{R}$ such that 
\begin{align}
\begin{bmatrix}
C^T & -1\\ -C_1^T &1\\ \vdots \\ -C_m^T &1 
\end{bmatrix} \cdot \begin{bmatrix}d\\ d_H \end{bmatrix} \prec \begin{bmatrix}0\\0\\ \vdots \\ 0 \end{bmatrix}
\end{align} 
\end{remark}

In the following, let us assume w.l.o.g that $C = 0_{n \times 1}$ and that $C \in \text{conv}\{C_1, \hdots, C_m\}$. We have:
\begin{align}\label{E3.10}
(h-f)(x) & = -\|x\|^2 + \max_{k \in \{1, \hdots, m\}} \|x - C_k\|^2 - r_k^2 \nonumber \\ 
&= \max_{k \in \{1, \hdots, m\}} \hspace{0.2cm} -2\cdot C_k^T \cdot X +\|C_k\|^2 - r_k^2 
\end{align} 
Let 
\begin{align}
x_1^{\star} \in \mathcal{X}^{\star}_1 = \mathop{\text{argmin}}_{h(x) \leq 1} \hspace{0.2cm} (h-f)(x) 
\end{align} 

There are three possibilities now:
\begin{enumerate}
\item $\mathcal{X}_1^{\star} \cap \{x | h(x) < 0\} = \emptyset$ 
\item $\mathcal{X}_1^{\star} \subset \{ x | h(x) < 0\}$ 
\item $ \exists x_1^{\star} \in \mathcal{X}_1^{\star} \cap \{x | h(x) = 0\}$
\end{enumerate}
\subsubsection{The case $\mathcal{X}_1^{\star} \cap \{x | h(x) < 0\} = \emptyset$  }  In this situation we can apply the Theorem \ref{T2.1} to solve (\ref{E3.1})
\subsubsection{The case $\mathcal{X}_1^{\star} \subset \{ x | h(x) < 0\}$} \label{ss3.2.2}
Note that (\ref{E3.10}) can be used to define a family of polytopes. Indeed, for some $R \in \mathbb{R}^{n \times 1}$, $(h-f)(x) \leq -R^2$ means
\begin{align}
\begin{bmatrix} -2\cdot C_1^T\\ \vdots\\ -2\cdot C_m^T \end{bmatrix} \cdot X + \begin{bmatrix} \|C_1\|^2-r_k^2 + R^2\\ \vdots\\ \|C_m\|^2 - r_m^2 + R^2\end{bmatrix} \preceq 0_{m \times 1}
\end{align} Let 
\begin{align}
\mathcal{P}_{R^2} = \left\{ x \in \mathbb{R}^{n \times 1} \biggr| (h-f)(x) \leq -R^2 \right\}
\end{align} 

We have the following remark about the polytope $\mathcal{P}_{R^2}$:
\begin{remark}\label{R3.2} 
Please note that for $R = 0$ we have that $\{x | h(x) \leq 0\} \subseteq \mathcal{P}_0$. Indeed, let $x \in \{x | h(x) \leq 0\}$ then $max_{k} \|x - C_k\|^2 - r_k^2 \leq 0$. On the other hand, $x \in \mathcal{P}_0$ if $(h-f)(x) \leq 0$. From (\ref{E3.10}) one can see that this easily happens. 

Furthermore, because $h(x_1^{\star}) < 0$ for all $x^{\star}_1 \in \mathcal{X}^{\star}_1$, follows that $h-f$ is bounded below, and we can say that $\exists \bar{R} \in \mathbb{R}$, the smallest, such that $\mathcal{P}_{R^2} = \emptyset$ for all $R > \bar{R}$. 

\textit{Increasing $R$, the set $\mathcal{P}_{R^2}$ shrinks from initially including the set $\{x | h(x) \leq 0\}$, for $R = 0$, to some points strictly inside $\{x| h(x) \leq 0\}$, for $R = \bar{R}$, then vanishes.}   
\end{remark}

In the following, we are interested in the largest $R$ such that $\mathcal{P}_{R^2} \cap \{x | h(x) = 0\} \neq \emptyset$.
 We shall call this $R^{\star}$ and give the following theorem:

\begin{theorem}  \label{T3.3} 
\begin{align} \label{E3.14}
R^{\star} = \max \hspace{0.2cm} \|x\|^2 \hspace{0.2cm} \text{s.t} \hspace{0.2cm} x \in \bigcap_{k = 1}^m \bar{\mathcal{B}}(C_k, r_k) 
\end{align}
\end{theorem}
\begin{proof} The problem (\ref{E3.1}) is maximizing a convex function over a convex domain, hence the maximizer will be on the boundary. For some $0 < R < R^{\star}$  we shall prove that $ \left( \{x | h(x) = 0\} \setminus \bar{\mathcal{B}}(0_{n \times 1},R) \right) \subseteq \mathcal{P}_{R^2}$, which means that the points on the boundary of $\bigcap_{k} \mathcal{B}(C_k,r_k)$ of magnitude greater than $R$, are inside $\mathcal{P}_{R^2}$. Indeed, this is easy to verify: let $h(x) = 0$ and $\|x \| > R$ then 
\begin{align}
(h-f)(x) \leq -R^2 \iff h(x) + R^2 - \|x\|^2 \leq 0 \iff x \in \mathcal{P}_{R^2}
\end{align} It follows that $\not\exists x_2$ with $h(x_2) = 0$, i.e, on the boundary of $\bigcap_k \bar{\mathcal{B}}(C_k,r_k)$, and $\|x_2\| > R^{\star}$, otherwise $\|x_2\| = R_2 > R^{\star}$ and $x_2 \in \{x|h(x) = 0\} \cap \mathcal{P}_{R_2^2} \neq \emptyset$ which is a contradiction with $R^{\star}$ being the largest with this property. 
\end{proof}

As a final remark, we say that finding $R^{\star}$ might be difficult in general. However, this can be done to a certain precision for \textbf{Type 3} problems, see below.
\subsubsection{The case $ \exists x_1^{\star} \in \mathcal{X}_1^{\star} \cap \{x | h(x) = 0\}$} \label{ss3.2.3}In this situation we argue that $R^{\star} = \bar{R}$ where $\bar{R}$ is given in Remark \ref{R3.2} and $R^{\star}$ is given by (\ref{E3.14}). Indeed, similar to the proof of Theorem \ref{T3.3} follows that $\not\exists R > \bar{R}$ such that $\{x | h(x) = 0\}$ and $\|x\| = R$. Assuming otherwise leads to the fact that $\exists x_3 \in \{x | h(x) = 0\} \cap \mathcal{P}_{R^2} \neq \emptyset$. This is a contradiction simply because for all $R > \bar{R}$ one has $\mathcal{P}_{R^2} = \emptyset$. 

\subsection{Type 3}
Let us consider the real subset sum problem: for a given $S \in \mathbb{R}^{n \times 1}$ we ask if exists $x \in \{0,1\}^{n \times 1}$ such that $S^T \cdot x = 0$. 

 We formulate the well known optimization problem for some $\beta > 0$:
\begin{align}\label{E3.16}
\max &\hspace{0.2cm} \sum_{k=1}^n x_k \cdot (x_k - 1) + \beta \cdot \sum_{k=1}^n x_k \cdot s_k \nonumber \\
& \text{s.t} \hspace{0.2cm} \begin{cases} S^T \cdot x \leq 0\\ 0 \leq x_k \leq 1 \\
1_{n \times 1}^T \cdot x - \frac{1}{2} \geq 0 \end{cases}
\end{align} where $x = \begin{bmatrix}x_1 &\hdots &x_n \end{bmatrix}^T \in \mathbb{R}^{n \times 1}$. Please note that due to the last constraint, the origin $0_{n \times 1}$ does not belong to the search space. 

As presented here \cite{sahni_comp} note that upon solving (\ref{E3.16}) the answer to the problem is positive if and only if the maximum is zero.  Indeed, on the search space, the objective function is less than or equal to zero and reaches its maximum value of zero on $x^{\star}$ if and only if $x^{\star} \in \{0,1\}^{n\times 1}$ and $S^T \cdot x^{\star} = 0$.  

Upon writing 
\begin{align}
& \sum_{k=1}^n x_k \cdot (x_k - 1) + \beta \cdot \sum_{k=1}^n x_k \cdot s_k =  x^T \cdot x  + (\beta \cdot S - 1_{n \times 1})^T \cdot x = \nonumber \\
& =  x^T \cdot x + 2 \cdot x^T \cdot \frac{\beta \cdot S - 1_{n \times 1}}{2} + \frac{1}{4} \|\beta \cdot S - 1_{n \times 1}\|^2 - \frac{1}{4} \|\beta \cdot S - 1_{n \times 1}\|^2 \nonumber \\
& = \left\| x - \frac{1}{2} (1_{n \times 1} - S) \right\|^2 - \frac{1}{4} \|\beta \cdot S - 1_{n \times 1}\|^2 
\end{align} it is obvious that we can solve (\ref{E3.16}) if we solve: 
\begin{align}\label{E3.18}
\max &\hspace{0.2cm} \left\| x - \frac{1}{2} \left(1_{n \times 1} - \frac{\beta}{\|S\|} \cdot S \right) \right\|^2 \nonumber \\
& \text{s.t} \hspace{0.2cm} \begin{cases} S^T \cdot x \leq 0\\ 0 \leq x_k \leq 1 \\ 
1_{n \times 1}^T \cdot x - \frac{1}{2} \geq 0\end{cases} 
\end{align}  In the following we focus on solving the problem (\ref{E3.18}). 

Let $\mathcal{P}$ denote the search space. This describes a polytope being the intersection of the unit hypercube with $\{x | S^T \cdot x \leq 0\}$ and $\left\{ x |  1_{n \times 1}^T \cdot x - \frac{1}{2} \geq 0 \right\}$ half spaces. Our first step is to approximate this polytope with a suitable intersection of n-disks. 
\subsubsection{Construction of a suitable intersection of n-disks}
For approximating the unit hypercube, we consider the following n-disks $\mathcal{B}(C_{k+}, r_{k+})$ and $\mathcal{B}(C_{k-}, r_{k-})$ where
\begin{align}
C_{k\pm} = \frac{1}{2} \cdot 1_{n \times 1} \pm q_k \cdot e_k
\end{align} where $e_k$ is the $k$'th column of the unit matrix in $\mathbb{R}^{n \times 1}$ and $q_k$ and $r_k$ are to be determined such that 
\begin{align}\label{E3.20}
\left\| \frac{1}{2} \cdot 1_{n \times 1} - \frac{1}{2} \cdot e_k - C_{k+} \right\|^2 + \left( \frac{\sqrt{n-1}}{2} \right)^2 = r_{k+}^2 \nonumber \\
\left\| \frac{1}{2} \cdot 1_{n \times 1} + \frac{1}{2} \cdot e_k - C_{k-} \right\|^2 + \left( \frac{\sqrt{n-1}}{2} \right)^2 = r_{k-}^2
\end{align} These conditions assure that each corner of the unit hypercube belongs to the boundary of $\bigcap_{k=1}^n  \left( \mathcal{B}(C_{k+}, r_{k+}) \cap \mathcal{B}(C_{k-}, r_{k-})\right)$. Indeed, w.l.o.g let us consider the corners of the unit hypercube which belong to the facet $x_1 = 1$. These will be in the form $\begin{bmatrix}1\\ y \end{bmatrix}$ where $y \in \{0,1\}^{n-1}$.  Then we have 
\begin{align}
\left\| \begin{bmatrix}1\\ y \end{bmatrix} - \left( \frac{1}{2}\cdot 1_{n \times 1} + \frac{1}{2} \cdot e_1 \right)\right\|^2 = \left\| \begin{bmatrix} 0\\ y - \frac{1}{2} \cdot 1_{n\times 1}\end{bmatrix}\right\|^2 = \left(\frac{\sqrt{n-1}}{2} \right)^2
\end{align} Finally, we show that $\begin{bmatrix} 1\\ y\end{bmatrix} \in \mathcal{B}(C_{k-},r_{k-})$. 
\begin{align}
\left\| \begin{bmatrix}1 \\ y \end{bmatrix} - C_{k-} \right\|^2 &= \left\|\begin{bmatrix}1\\ y \end{bmatrix} - \left( \frac{1}{2} \cdot 1_{n \times 1} + \frac{1}{2}\cdot e_1 \right) +  \left( \frac{1}{2} \cdot 1_{n \times 1} + \frac{1}{2}\cdot e_1 \right) - C_{k-}\right\|^2 \nonumber \\
& = \left\|\begin{bmatrix}0\\ y - \frac{1}{2} \cdot 1_{n \times 1} \end{bmatrix} \right\|^2 +\left\| \frac{1}{2} \cdot 1_{n \times 1} + \frac{1}{2} \cdot e_1 - C_{k-} \right\|^2  = r_{k-}^2
\end{align} because $ \frac{1}{2} \cdot 1_{n\times 1} + \frac{1}{2} \cdot e_1 - C_{k-} = e_1 \cdot \left( \frac{1}{2} + q_{k}\right) $ which is orthogonal on $\begin{bmatrix} 0\\ y - \frac{1}{2} \cdot 1_{n \times 1} \end{bmatrix}$. Therefore the unit hypercube is included in the intersection of n-disks and has the corners touching the boundary of the intersection. 

Next we approximate the half space $\{x | S^T\cdot x \leq 0\}$ with the closed n-disk $\bar{\mathcal{B}}(C_s, r_s)$ where we compute $C_s$ and $r_s$ as follows. Let $P_s = \frac{1}{2} \cdot 1_{n \times 1} + \frac{t}{\|S\|} \cdot S$ such that $S^T \cdot P_s = 0$ hence $S^T \cdot 1_{n \times 1} + \frac{t}{\|S\|} \cdot \|S\|^2 = 0$ and consequently $P_s = \frac{1}{2} \cdot 1_{n \times 1} - \frac{S^T \cdot 1_{n \times 1}}{\|S\|^2} \cdot S$. Let 
\begin{align}
C_s &= P_s - \frac{\tilde{q}_s}{\|S\|} \cdot S = \frac{1}{2} \cdot 1_{n \times 1} - \left( \frac{S^T \cdot 1_{n \times 1} }{\|S\|} + \tilde{q}_s\right)  \cdot \frac{S}{\|S\|} \nonumber \\
& = \frac{1}{2} \cdot 1_{n \times 1} - q_s \cdot \frac{S}{\|S\|}
\end{align} and 
\begin{align}\label{E3.24}
\tilde{r}_s^2 = \frac{n}{4} - \left\|\frac{1}{2} \cdot 1_{n \times 1} - P_s \right\|^2 = \frac{n}{4} - \left( \frac{S^T \cdot 1_{n \times 1}}{\|S\|}\right)^2
\end{align}
 and we require the following relation to hold:
\begin{align}\label{E3.25}
\| P_s - C_s\|^2 + \tilde{r}_s^2 = r_s^2
\end{align} This is need to ensure that if a corner of the unit hypercube belongs to the hyperplane $\{x | S^T \cdot x = 0\}$ then it will also belong to the hyper sphere $\{x | \|x - C_s\| = r_s\}$.  Indeed, let $x \in \{ 0,1\}^{n \times 1}$ such that $x^T\cdot S = 0$, then we show that $\|C_s - x\| = r_s$. 
\begin{align}
\|C_s - x\|^2 &= \left\|C_s - P_s + P_s - x \right\|^2 = \| C_s -P_s\|^2 + \|P_s - x\|^2 + 2 \cdot (C_s-P_s)^T\cdot (P_s-x) \nonumber \\
&= r_s^2 - \tilde{r}_s^2 + \left\| x - \frac{1}{2}\cdot 1_{n \times 1}\right\|^2 - \left\| \frac{1}{2}\cdot 1_{n \times 1} - P_s \right\|^2 - 2 \cdot \tilde{q}_s\cdot\frac{S^T}{\|S\|} \cdot (P_s - x) \nonumber \\
& = r_s^2 - \tilde{r}_s^2 + \frac{n}{4} -  \left\| \frac{1}{2}\cdot 1_{n \times 1} - P_s \right\|^2 = r_s^2
\end{align} because $S^T \cdot P_s = 0 = S^T \cdot x$ and $\|P_s - x\|^2 + \left \| \frac{1}{2} \cdot 1_{n \times 1} - P_s \right\|^2 = \left\| x - \frac{1}{2}\cdot 1_{n \times 1} \right\|^2$ by applying Pitagora's theorem in the right triangle $\triangle x \frac{1_{n \times 1}}{2} P_s$. Finally please note that $\left\|x - \frac{1_{n \times 1}}{2}\right\| = \frac{\sqrt{n}}{2}$ for all $x \in \{0,1\}^{n \times 1}$. 

The last constraint is the half space $\left\{x | 1_{n\times 1}^T \cdot x - \frac{1}{2} \geq 0 \right\}$ which will be "approximated" with the closed n-disk $\bar{\mathcal{B}}(C_h, r_h)$ with $C_h$ and $r_h$ to be calculated below. Let $P_h = \frac{1}{2} \cdot 1_{n \times 1} + t \cdot \frac{1_{n \times 1}}{\|1_{n \times 1}\|}$ such that $1_{n \times 1}^T \cdot P_h - \frac{1}{2} = 0$. It follows that $1_{n \times 1}^T \cdot \left( \frac{1}{2} \cdot 1_{n \times 1} + t \cdot \frac{1_{n \times 1}}{\|1_{n \times 1}\|} \right) - \frac{1}{2} = 0$ hence $P_h = \frac{1}{2} \cdot 1_{n \times 1} + \frac{1-n}{2\cdot \sqrt{n}} \cdot \frac{1_{n \times 1}}{\|1_{n \times 1}\|} = \frac{1}{2\cdot n} \cdot 1_{n \times 1} $ then

\begin{align}
C_h &= P_h + \tilde{q}_h\cdot \frac{1_{n \times 1}}{\|1_{n\times 1}\|} = \frac{1}{2}\cdot 1_{n \times 1} + \left( \frac{1-n}{2 \cdot \sqrt{n}} + \tilde{q}_h \right) \cdot \frac{1_{n\times 1}}{\|1_{n \times 1}\|} \nonumber \\
& = \frac{1}{2}\cdot 1_{n\times 1} + q_h \cdot \frac{1_{n\times 1}}{\|1_{n \times 1}\|}
\end{align} Let
\begin{align} \label{E3.28}
\tilde{r}_h^2 = \frac{n}{4} - \left\| \frac{1}{2} \cdot 1_{n \times 1} - P_h \right\|^2 = \frac{n}{4} - \left( \frac{1}{2} - \frac{1}{2\cdot n}\right)^2 \cdot n = \hdots = \frac{1}{2} - \frac{1}{4\cdot n}
\end{align}

As above, we require $q_h$ and $r_h$ to meet the following constraint:
\begin{align}\label{E3.29}
\|C_h - P_h\|^2 + \tilde{r}_h^2  = r_h^2
\end{align} This is enough to ensure that $\left\{x | \left\|x - \frac{1}{2} \cdot 1_{n\times 1} \right\| \leq \frac{\sqrt{n}}{2}, 1_{n\times 1}^T \cdot x - \frac{1}{2} \geq 0 \right\} \subseteq \bar{\mathcal{B}}(C_h, r_h)$ as we will see in the next subsection.

Finally, we choose 
\begin{align}\label{E3.30a}
q_{k \pm} = q_h = q_s = \rho
\end{align}  and define: 

\begin{align}
\mathcal{C}_{\rho} = \bigcap_{k=1}^n \left( \bar{\mathcal{B}}(C_{k+}, r_{k+}) \cap \bar{\mathcal{B}}(C_{k-}, r_{k-})\right) \cap \bar{\mathcal{B}}(C_s, r_s) \cap \bar{\mathcal{B}}(C_h, r_h)
\end{align} 

Please note that there is a one to one relation between $\rho$ and the radii of the n-disks. If the radii are constrained to be greater than some fixed value, then this is achievable by choosing a sufficiently large $\rho$. 

\subsubsection{Analysis of the set $\mathcal{C}_{\rho}$}

We begin this subsection with a very useful lemma:

\begin{lemma} \label{L3.4} Let us consider the following n-disks 
\begin{align}
\mathcal{D}_1 &= \{x \in \mathbb{R}^{n\times 1} | \|-q_1\cdot e_1 - x\| \leq r_1\} \nonumber \\
\mathcal{D}_2 &= \{x \in \mathbb{R}^{n\times 1}| \|-q_2 \cdot e_1 - x\| \leq r_2\} \nonumber \\
\mathcal{D}_3 &= \{x \in \mathbb{R}^{n\times 1} | \|-q_3 \cdot e_1 - x\| \leq r_3\}
\end{align} with $q_1, q_2 > 0$ and $r_1 > r_2 > r_3 \geq 0$ such that exists 
\begin{align}
0<a^2 = r_1^2 - q_1^2  = r_2^2 - q_2^2 = r_3^2 - q_3^2
\end{align} that is the n-disks share a common n-1 sphere i.e $\{x | e_1^T\cdot x = 0, \|x\| = a\}$. Let $\mathcal{H} = \{x \in \mathbb{R}^{n\times 1} | e_1^T \cdot x \leq 0\}$ and $\mathcal{G} = \{x \in \mathbb{R}^{n\times 1} | e_1^T \cdot x \geq 0\} = \mathbb{R}^{n\times 1} \setminus \text{int}(H)$. Then the following inclusions are true
\begin{enumerate}
\item
\begin{align}
\mathcal{H} \cap \mathcal{D}_3 \subseteq \mathcal{D}_1 \cap \mathcal{D}_3 \subseteq \mathcal{D}_2 \cap \mathcal{D}_3 \subseteq \mathcal{D}_2 \nonumber 
\end{align}
\item 
\begin{align}
\mathcal{G} \cap \mathcal{D}_1 \subseteq \mathcal{G} \cap \mathcal{D}_2 \subseteq \mathcal{G} \cap \mathcal{D}_3
\end{align}
\end{enumerate}
\end{lemma}
\begin{proof} see Appendix 
\end{proof}

\begin{remark}  \label{R3.5}
Using the above Lemma \ref{L3.4} and the construction of the set $\mathcal{C}_{\rho}$ the following can be proven for $\rho = \min\{q_{k\pm}, q_s, q_h\}$
\begin{enumerate}
\item Exists $\bar{\rho} > 0$ such that $\min \{r_{k\pm}, r_s, r_h\} \geq \frac{\sqrt{n}}{2}$, $C_s^T \cdot S < 0$ and $C_h^T \cdot 1_{n\times 1} > 0$ for all $\rho \geq \bar{\rho}$. This is easy to prove using the construction of the n-disks. 
\item If $\rho > \bar{\rho}$ then one has
\begin{align}
 \mathcal{P} \subseteq \mathcal{C}_{\rho} 
\end{align} 
 Indeed, since we know that $\mathcal{P} \subseteq \bar{\mathcal{B}}\left( \frac{1}{2} \cdot 1_{n \times 1}, \frac{\sqrt{n}}{2}\right) $ in order to prove that $\mathcal{P} \subseteq \mathcal{C}_{\rho}$ is enough to prove that any half space composing $\mathcal{P}$ intersected with $\bar{\mathcal{B}}\left( \frac{1}{2} \cdot 1_{n \times 1}, \frac{\sqrt{n}}{2}\right)$ is included in an n-disk composing $\mathcal{C}_{\rho}$ intersected with $\bar{\mathcal{B}}\left( \frac{1}{2} \cdot 1_{n \times 1}, \frac{\sqrt{n}}{2}\right)$. W.l.o.g take in the above Lemma \ref{L3.4} $\mathcal{H} = \{ x | S^T \cdot x \leq 0\}$, $\mathcal{D}_3 = \bar{\mathcal{B}}\left( \frac{1}{2} \cdot 1_{n \times 1}, \frac{\sqrt{n}}{2}\right)$ and $\mathcal{D}_1 = \bar{\mathcal{B}}(C_s, r_s)$ to obtain $\mathcal{H} \cap \mathcal{D}_3 \subseteq \mathcal{D}_1 \cap \mathcal{D}_3$
\item For $\rho > \bar{\rho}$ one also has 
\begin{align}
\mathcal{C}_{\rho} \subseteq \bar{\mathcal{B}}\left( \frac{1}{2} \cdot 1_{n \times 1}, \frac{\sqrt{n}}{2}\right) \hspace{0.5cm} 
\end{align} 
 Indeed
\begin{align}\label{E3.37a}
\mathcal{C}_{\rho} \subseteq \bigcap_{k = 1}^n \left( \bar{\mathcal{B}}(C_{k+}, r_{k+}) \cap  \bar{\mathcal{B}}(C_{k-}, r_{k-})\right) =\tilde{\mathcal{C}}_{\rho}
\end{align} W.l.o.g let us analyze
\begin{align}
\bar{\mathcal{B}}(C_{k+}, r_{k+}) \subseteq \left\{ x| e_k^T\cdot x \geq 0 \right\} \cup \left( \{ x | e_k^T\cdot x \leq 0\} \cap \bar{\mathcal{B}}(C_{k+}, r_{k+}) \right)
\end{align} It can be proven using the above Lemma \ref{L3.4} and the constructions of the n-disks $\bar{\mathcal{B}}(C_{k+}, r_{k+})$ that 
\begin{align}
\left( \{ x | e_k^T\cdot x \leq 0\} \cap \bar{\mathcal{B}}(C_{k+}, r_{k+}) \right) &\subseteq \left( \{ x | e_k^T\cdot x \leq 0\} \cap \bar{\mathcal{B}}\left(\frac{1}{2}\cdot 1_{n\times 1}, \frac{\sqrt{n}}{2} \right) \right) \nonumber \\
&\subseteq \bar{\mathcal{B}}\left(\frac{1}{2}\cdot 1_{n\times 1}, \frac{\sqrt{n}}{2} \right) 
\end{align} hence we obtained 
\begin{align}
\bar{\mathcal{B}}(C_{k+}, r_{k+}) \subseteq \left\{ x| e_k^T\cdot x \geq 0 \right\} \cup \bar{\mathcal{B}}\left(\frac{1}{2}\cdot 1_{n\times 1}, \frac{\sqrt{n}}{2} \right)
\end{align}

 Finally considering the following property regarding a finite intersection of sets for $A_k, B \subseteq \mathbb{R}^{n \times 1}$
\begin{align} 
\bigcap_{k=1}^n A_k \subseteq B \Rightarrow \bigcap_{k=1}^n (A_k \cup B) \subseteq B
\end{align}

Therefore, since the unit hypercube $\bigcap_{k=1}^n (\{x | e_k^T\cdot x \geq 0\} \cap \{ x | e_k^T\cdot x \leq 1\} ) \subseteq \bar{\mathcal{B}}\left(\frac{1}{2}\cdot 1_{n\times 1}, \frac{\sqrt{n}}{2} \right)$ one obtains $\tilde{\mathcal{C}} \subseteq \bar{\mathcal{B}}\left(\frac{1}{2}\cdot 1_{n\times 1}, \frac{\sqrt{n}}{2} \right)$. 
\end{enumerate}

\end{remark}

The following lemma can be proved
\begin{lemma}\label{L3.6}
Given $\delta > 0$ fixed, if $\exists \mathcal{B}(x,\delta) \subseteq \mathcal{P} $ then exists $\rho_{\delta} > 0$ such that for all $\rho \geq \rho_{\delta}$ one has 
\begin{align}
\mathcal{P} \subseteq \mathcal{C}_{\rho} \subseteq \bigcup_{x \in \mathcal{P}} \mathcal{B}(x,\delta)
\end{align}
\end{lemma}
\begin{proof} See Appendix
\end{proof}

What the above lemma says is that by using large enough radii for the n-disks involved (please note that the construction presented above allows the radii of the n-disks to be chosen arbitrarily large) one can assure a fixed desired level of approximation of $\mathcal{P}$. 

Let us define
\begin{align} \label{E3.42}
C = \frac{1}{2} \cdot 1_{n\times 1} - \frac{\beta}{2} \cdot \frac{S}{\|S\|}
\end{align} and consider the problem
\begin{align}
\max \hspace{0.2cm} \|x- C\| \hspace{0.2cm} \text{s.t} \hspace{0.2cm} x \in \mathcal{C}_{\rho}
\end{align}  for some $\rho > 0$. 
We continue the subsection with a very important lemma:
\begin{lemma}\label{L3.7}
For $\bar{\rho} < \frac{\beta}{2} < q_s = \rho$  let
\begin{align}
x^{\star} \in \mathcal{U}^{\star} = \mathop{\text{argmax}}_{x \in \mathcal{C}_{\rho}} \hspace{0.2cm} \|x - C\|^2
\end{align} then 
\begin{enumerate}
\item $\left\|x^{\star} - \frac{1}{2} \cdot 1_{n \times 1} \right\| \leq \frac{\sqrt{n}}{2}$
\item $\left\|x^{\star} - \frac{1}{2} \cdot 1_{n \times 1} \right\| = \frac{\sqrt{n}}{2} $ iff $x^{\star} \in \{0,1\}^{n\times 1}$
\item If exists $x_1 \in \{0,1\}^{n \times 1}$ with $S^T \cdot x_1 = 0$ then $x_1 \in \mathcal{U}^{\star}$
\item If exists $x_1 \in \{0,1\}^{n \times 1}$ with $S^T \cdot x_1 = 0$ then for all $x^{\star} \in \mathcal{U}^{\star}$ one has $x^{\star} \in \{0,1\}^{n \times 1}$ and $S^T \cdot x^{\star} = 0$
\end{enumerate} 
%
\end{lemma}
\begin{proof} Since $\rho > \bar{\rho}$, from Remark \ref{R3.5} follows $\mathcal{C}_{\rho} \subseteq \bar{\mathcal{B}}\left(\frac{1}{2}\cdot 1_{n\times 1}, \frac{\sqrt{n}}{2}\right)$ Next: 
\begin{enumerate}
\item The first claim follows easily from the fact that $x^{\star} \in \mathcal{C}_{\rho} \subseteq \bar{\mathcal{B}}\left(\frac{1}{2}\cdot 1_{n\times 1}, \frac{\sqrt{n}}{2}\right)$ 
\item For the second claim it is easy to verify the reverse implication. We focus on the direct one. Let $\left\|x^{\star} - \frac{1}{2} \cdot 1_{n \times 1} \right\| = \frac{\sqrt{n}}{2} $. Since 
\begin{align}
x^{\star} \in \mathcal{C}_{\rho} \cap \partial \bar{\mathcal{B}}\left(\frac{1}{2}\cdot 1_{n\times 1}, \frac{\sqrt{n}}{2}\right) \subseteq \tilde{C}_{\rho} \cap \bar{\mathcal{B}}\left(\frac{1}{2}\cdot 1_{n\times 1}, \frac{\sqrt{n}}{2}\right)
\end{align} follows that  $x^{\star} \in \{0,1\}^{n \times 1}$ since it can be proven that for $\rho > \bar{\rho}$ one has $\tilde{C}_{\rho} \cap \bar{\mathcal{B}}\left(\frac{1}{2}\cdot 1_{n\times 1}, \frac{\sqrt{n}}{2}\right) \subseteq \{0,1\}^{n \times 1}$. See (\ref{E3.37a}) for the definition of $\tilde{C}_{\rho}$. It is basically an approximation of the unit hypercube. 
\item Let $x_1 \in \{0,1\}^{n\times 1} $ with $S^T\cdot x_1 = 0$ then by the construction we get $x_1 \in \mathcal{C}_{\rho}$ then $\left\|x_1 - \frac{1}{2} \cdot 1_{n \times 1} \right\| = \frac{\sqrt{n}}{2}$ and $\|C_s - x_1\| = r_s$. We shall prove that $\|x^{\star} - C\| \leq \|x_1 - C\|$. 
%

Since $q_s > \frac{\beta}{2} > \bar{\rho}$ we can choose in Lemma \ref{L3.4} $\mathcal{D}_1 = \bar{\mathcal{B}}(C_s, r_s)$, $\mathcal{D}_3 = \bar{\mathcal{B}}\left( \frac{1}{2} \cdot 1_{n \times 1}, \frac{\sqrt{n}}{2}\right)$ and $\mathcal{D}_2 = \bar{\mathcal{B}}(C, \|C - x_1\|)$ to  assert that $\mathcal{D}_1 \cap \mathcal{D}_3 \subseteq \mathcal{D}_2$. That is, for all $x \in \mathcal{C}_{\rho} \subseteq \mathcal{D}_3 \cap \mathcal{D}_1$ one has $x \in \mathcal{D}_2$ hence $\|x^{\star}-C\| \leq \|C - x_1\|$. 
\item It is known from the previous point that $ \|C - x^{\star}\| = \|C - x_1\|$. Then it follows that $\|C_s - x^{\star}\| = r_s$ and $\left\| \frac{1}{2}\cdot 1_{n\times 1} - x^{\star} \right\| = \frac{\sqrt{n}}{2}$. Indeed, assume that $\|C_s - x^{\star}\| < r_s$ and/or $\left\| x - \frac{1}{2} \cdot 1_{n \times 1}\right\| < \frac{\sqrt{n}}{2}$. Then by the use of Lemma \ref{L3.4}, similarly to previous point, one obtains $\|C - x^{\star}\| < \|C - x_1\|$, which is false. Using both is enough to prove that $S^T\cdot x^{\star} = 0$ and using the later is enough to prove, using the second point in this enumerate, that $x^{\star} \in \{0,1\}^{n \times 1}$
\end{enumerate}
\end{proof}

\subsubsection{Solution to some real subset sum problems}
Let us consider the problem:
\begin{align}
\mathcal{V}^{\star} = \mathop{\text{argmax}}_{x \in \mathcal{P}} \hspace{0.2cm} \|x - C\|^2
\end{align}

For $\bar{\rho} < \frac{\beta}{2} < q_s = \rho$  let
\begin{align} \label{E3.47}
\mathcal{U}^{\star} = \mathop{\text{argmax}}_{x \in \mathcal{C}_{\rho}} \hspace{0.2cm} \|x - C\|^2
\end{align} 
We reason as follows: if $\mathcal{V}^{\star} \subseteq \left( \{0,1\}^{n \times 1} \cap \{x | S^T \cdot x =0\} \right) $ then from Lemma \ref{L3.7} point 4, follows immediately that $\mathcal{V}^{\star} = \mathcal{U}^{\star}$ and every point in $\mathcal{U}^{\star}$ has that property. In the following we show that if a point has this property we will find it. 


In order to solve (\ref{E3.47}) as described in Theorem \ref{T3.3} we form the function 
\begin{align}
&h(x) = \max_{ p \in \{k\pm, s, h\}} \|x - C_p\|^2 - r_p^2 
\end{align} 
where obviously $\mathcal{C}_{\rho} = \{x | h(x) \leq 0\}$. 
Let us define
\begin{align}
\mathcal{W}^{\star} = \mathop{\text{argmin}}_{h(x) \leq 1} \hspace{0.2cm} h(x) - \|x - C\|^2
\end{align} \textbf{ For the ease of presentation, in the following we assume that $\mathcal{W}^{\star}$ has only one element }, therefore $\mathcal{W}^{\star} = \{w^{\star}\}$. If $h(w^{\star}) > 0$ then we can apply Theorem \ref{T2.1} to solve (\ref{E3.47}). Otherwise if $h(w^{\star}) = 0$  then using Section \ref{ss3.2.3} one can eventually prove that $w^{\star} \in \mathcal{U}^{\star}$. Finally, the most exciting case (and likely to meet) is when $h(w^{\star}) < 0$.  This is what we treat in the following.
\begin{remark} The assumption that $\mathcal{W}^{\star} = \{ w^{\star}\}$ actually limits the choice of $S$, therefore this algorithm might not be applicable for any $S \in \mathbb{R}^{n \times 1}$. 
\end{remark}

 So the case $h(w^{\star}) < 0$ is to be analyzed. For start we assume something slightly stronger, i.e that $\exists \mathcal{B}(w^{\star}, \epsilon>0) \subseteq \mathcal{P} \subseteq \mathcal{C}_{\rho}$. 

In this case, obviously $h(w^{\star}) <0 $ and we apply the theory presented in Section \ref{ss3.2.2}. Let us form the family of polytopes
\begin{align}
\mathcal{P}_{R^2} = \{x | h(x) - \|x - C\|^2 \leq -R^2\}
\end{align} for any $R \geq 0$. According to Section \ref{ss3.2.2} and Theorem \ref{T3.3} in order to solve (\ref{E3.47}) we should increase $R$ until $\mathcal{P}_{R^2} \subseteq \mathcal{C}_{\rho}$. The smallest $R$ for which this happens shall be denoted $R^{\star}$ and $\mathcal{U}^{\star} = \partial \mathcal{P}_{(R^{\star})^2} \cap \partial \mathcal{C}_{\rho}$.

The difficulty here is, off course, asserting if a polytope (a.k.a $\mathcal{P}_{R^2}$) is included in an intersection of n-disks (a.k.a $\mathcal{C}_{\rho}$). Actually this problem seems as, if not even more, complicated than the initial one. However, in this special case, i.e when $\partial \mathcal{P} \cap \partial \mathcal{C}_{\rho} \neq \emptyset$ and other properties that $\mathcal{P}$ and $\mathcal{C}_{\rho}$ share, we show that this can actually be done in polynomial time.

\textit{The obvious observation here is the following: it is easy to assert if $\mathcal{P}_{R^2} \subseteq \mathcal{P}$!}. 

We already know that $\mathcal{P} \subseteq \mathcal{C}_{\rho}$, so having $\mathcal{P}_{R^2} \subseteq \mathcal{P}$ is sufficient to say that $\mathcal{P}_{R^2} \subseteq \mathcal{C}_{\rho}$. We formulate the following question:

\textit{Is $R^{\star}$ the smallest $R$ for which $\mathcal{P}_{R^2} \subseteq \mathcal{P}$ ? }

Since $\mathcal{P} \subseteq \mathcal{C}_{\rho} \subseteq \mathcal{P}_{0}$ and $w^{\star} \in \mathcal{P}$ follows that exists $\tilde{R}^{\star}$ such that $\tilde{R}^{\star}$ is the smallest for which $\mathcal{P}_{R^2} \subseteq \mathcal{P}$. It is obvious that $\tilde{R}^{\star} \geq R^{\star}$. The next question is:

\textit{For $R = R^{\star}$ exists $x \in \mathcal{P}_{R^2} \setminus \mathcal{P}$ ? }

If the answer to the above question is $NO$ then it follows that $\tilde{R}^{\star} = R^{\star}$. We shall prove that this is the case. Let us assume that exists 
\begin{align} \label{E3.52b}
x_0 \in \mathcal{P}_{(R^{\star})^2} \setminus \mathcal{P}
\end{align} and exists $\delta_0 > 0$ such that $\mathcal{B}(x_0, \delta_0) \subseteq \mathcal{P}_{(R^{\star})^2} \setminus \mathcal{P}$.

In this case, starting from $\mathcal{C}_{\rho}$ we form another intersection of  n-disks as follows:
\begin{align}
\hat{\mathcal{C}}_{\hat{\rho}} = \bigcap_{k=1}^n \left( \bar{\mathcal{B}}(\hat{C}_{k+}, \hat{r}_{k+}) \cap \bar{\mathcal{B}}(\hat{C}_{k-}, \hat{r}_{k-})\right) \cap \bar{\mathcal{B}}(\hat{C}_s, \hat{r}_s) \cap \bar{\mathcal{B}}(\hat{C}_h, \hat{r}_h)
\end{align} where for $\alpha > 0$ we have:
\begin{align} \label{E3.52}
\hat{C}_{k\pm} &= \frac{1}{2}\cdot 1_{n \times 1} \pm \alpha \cdot q_k \cdot e_k \nonumber \\
 \hat{C}_s &= \frac{1}{2} \cdot 1_{n \times 1} - \alpha \cdot q_s \cdot \frac{S}{\|S\|} \nonumber \\
\hat{C}_h & = \frac{1}{2} \cdot 1_{n \times 1} + \alpha \cdot q_h \cdot \frac{1_{n \times 1}}{\|1_{n\times 1}\|}
\end{align} with the constraints on the radii 
\begin{align} \label{E3.53}
&\left\| \frac{1}{2} \cdot 1_{n \times 1} - \frac{1}{2} \cdot e_k - \hat{C}_{k+} \right\|^2 + \left( \frac{\sqrt{n-1}}{2}\right)^2 = \hat{r}_{k+}^2 \nonumber \\
& \left\| \frac{1}{2} \cdot 1_{n \times 1} + \frac{1}{2} \cdot e_k - \hat{C}_{k-} \right\|^2 + \left( \frac{\sqrt{n-1}}{2}\right)^2 = \hat{r}_{k-}^2 \nonumber \\
&\| \hat{C}_s - P_s\|^2 + \tilde{r}_s^2 = \hat{r}_s^2 \nonumber \\
&\| \hat{C}_h - P_h\|^2 + \tilde{r}_h^2 = \hat{r}_h^2
\end{align}

Given the fact that $\exists \mathcal{B}(w^{\star}, \epsilon) \subseteq \mathcal{P}$ we can apply Lemma \ref{L3.6} with $\delta = \min \left\{ \epsilon, \frac{\delta_0}{2} \right\}$ to assure the existence of $\rho_{\delta}$ such that if $\hat{\rho}  \geq \rho_{\delta}$ then 
\begin{align} \label{E3.54}
\mathcal{P} \subseteq \hat{\mathcal{C}}_{\hat{\rho}} \subseteq \bigcup_{x \in \mathcal{P}} \mathcal{B}\left(x, \frac{\delta_0}{2} \right)
\end{align} Please note that we can achieve $\hat{\rho} = \alpha \cdot \rho \geq \rho_{\delta}$ by increasing $\alpha$ in (\ref{E3.52} and \ref{E3.53}).
\begin{remark} \label{R3.9b}
From the equation (\ref{E3.54}) one has $x_0 \not\in \hat{\mathcal{C}}_{\hat{\rho}}$
\end{remark}

Finally we define inspired from (\ref{E3.42})
\begin{align}
\hat{C} = \frac{1}{2} \cdot 1_{n \times 1} - \alpha \cdot \frac{\beta}{2} \cdot \frac{S}{\|S\|}
\end{align} and consider the problem:
\begin{align}
\hat{\mathcal{U}}^{\star} = \mathop{\text{argmax}}_{x \in \hat{C}_{\hat{\rho}}} \hspace{0.2cm} \|x - \hat{C}\|^2
\end{align} In order to solve this problem, as above, we define:
\begin{align}
\hat{h}(x) = \max_{p \in \{k\pm, s, h\}} \hspace{0.2cm} \|x - \hat{C}_p\|^2 - \hat{r}_p^2
\end{align} with the obvious remark that $\hat{\mathcal{C}}_{\hat{\rho}} = \{x | \hat{h}(x) \leq 0\}$. Then finally define
\begin{align}
\hat{\mathcal{W}}^{\star} = \mathop{\text{argmin}}_{\hat{h}(x) \leq 1} \hat{h}(x) - \|x - \hat{C}\|^2
\end{align} 

The last main contribution of this paper is the following lemma:
\begin{lemma}\label{L3.10}
For $R \geq R^{\star}$ exists $\hat{R}$ such that 
\begin{align}
\mathcal{P}_{R^2} =\hat{\mathcal{P}}_{\hat{R}^2}
\end{align} where
\begin{align}
\mathcal{P}_{R^2} &= \{x | h(x) - \|x - C\|^2 \leq -R^2\}  \nonumber \\ \hat{\mathcal{P}}_{\hat{R}^2} &= \{x | \hat{h}(x) - \|x - \hat{C}\|^2 \leq - \hat{R}^2\}
\end{align}
\end{lemma}
\begin{proof} 
See proof for Lemma \ref{LB.1} in Appendix. 
\end{proof}

Using the above Lemma \ref{L3.10} it is easy to prove that $\hat{\mathcal{W}}^{\star} \in \text{int}(\mathcal{P}) \subseteq \hat{\mathcal{C}}_{\hat{\rho}}$ hence we can define 

\begin{align}
\hat{R}^{\star} = \inf \{\hat{R} | \mathcal{P}_{\hat{R}^2} \subseteq \hat{\mathcal{C}}_{\hat{\rho}}\}
\end{align}
\begin{corollary}
 If  $ \exists x \in \mathcal{P} \cap \{0,1\}^{n\times 1} \cap \{x | S^T \cdot x = 0\}$ then 
\begin{align}
\mathcal{P}_{(R^{\star})^2} = \hat{\mathcal{P}}_{\left(\hat{R}^{\star}\right)^2}
\end{align}
\end{corollary}
\begin{proof}
If  $ \exists x \in \mathcal{P} \cap \{0,1\}^{n\times 1} \cap \{x | S^T \cdot x = 0\}$ then 
$ \exists x \in \hat{\mathcal{C}}_{\hat{\rho}} \cap \{0,1\}^{n\times 1} \cap \{x | S^T \cdot x = 0\} $ and follows that $\hat{\mathcal{U}}^{\star} \subseteq \{0,1\}^{n\times 1} \cap \{x | S^T \cdot x = 0\} $. It is also known that $\forall x \in \hat{\mathcal{U}}^{\star} $ one has $ x \in \partial \hat{\mathcal{P}}_{\left( \hat{R}^{\star} \right)^2} $

Since for all $\hat{R}_1 < \hat{R}_2$ one has $\hat{\mathcal{P}}_{R_2}^2 \cap \partial \hat{\mathcal{P}}_{R_1^2} = \emptyset$ it is easy to prove that $\mathcal{P}_{\left(R^{\star} \right)^2} = \mathcal{P}_{\left(\hat{R}^{\star} \right)^2}$. Indeed, for $R = R^{\star}$ let $\hat{R}_1$ be given by the Lemma \ref{L3.10} such that 
\begin{align}
\mathcal{P}_{(R^{\star})^2} = \hat{\mathcal{P}}_{\left(\hat{R}_1\right)^2}
\end{align} It is known that $ \forall x \in \{0,1\}^{n \times 1} \cap \{x | S^T \cdot x = 0\} \cap \mathcal{P}$ one has $x \in\partial \mathcal{P}_{(R^{\star})^2} = \partial \hat{\mathcal{P}}_{\hat{R}_1^2}$. Assuming $\hat{R}_1 < \hat{R}^{\star}$ follows $\hat{\mathcal{P}}_{(\hat{R}^{\star})^2} \cap \partial \hat{\mathcal{P}}_{(R_1)^2} = \emptyset$ hence $x \not\in \hat{\mathcal{P}}_{(\hat{R}^{\star})^2}$ which is a contradiction. On the other hand, assuming $\hat{R}_1 > \hat{R}^{\star}$ follows that $\hat{\mathcal{P}}_{(\hat{R}_1)^2} \cap \partial \hat{\mathcal{P}}_{(\hat{R}^{\star})^2} = \emptyset$ hence $x \not\in \hat{\mathcal{P}}_{(\hat{R}^{\star})^2}$, which is again a contradiction. It follows that $\hat{R}_1 = \hat{R}^{\star}$
\end{proof}.

Finally we conclude that $x_0$ given by (\ref{E3.52b}) cannot exist in $\mathcal{P}_{(R^{\star})^2} \subseteq \hat{\mathcal{C}}_{\hat{\rho}}$ since according to Remark \ref{R3.9b} $x_0 \not\in \hat{\mathcal{C}}_{\hat{\rho}}$.  

\section{Acknowledgments}
I would like to thank prof. B. Gavrea for the helpful suggestions and corrections regarding the development of the theory for first framework and  the applications of the first and second framework for Type 1 and 2 problems. 

\section{Conclusion and future work}

We presented two theoretical frameworks which can be used to solve certain instances of maximization problems over strongly convex domains. 
\bibliographystyle{siamplain}

\appendix
\section{Intersection of large and small spheres}

Let us first give a more fundamental lemma:
\begin{lemma} \label{LA1} Let us consider the following n-disks 
\begin{align}
\mathcal{D}_1 &= \{x \in \mathbb{R}^{n\times 1} | \|-q_1\cdot e_1 - x\| \leq r_1\} \nonumber \\
\mathcal{D}_2 &= \{x \in \mathbb{R}^{n\times 1}| \|-q_2 \cdot e_1 - x\| \leq r_2\} 
\end{align} with $q_1 > 0$ and $r_1 > r_2 \geq 0$ such that exists 
\begin{align}\label{EA.2}
0<a^2 = r_1^2 - q_1^2  = r_2^2 - q_2^2 
\end{align} that is the n-disks share a common $n-1$ sphere i.e $\{x | e_1^T\cdot x = 0, \|x\| = a\}$. Let $\mathcal{H} = \{x \in \mathbb{R}^{n\times 1} | e_1^T \cdot x \leq 0\}$ and $\mathcal{G} = \{x \in \mathbb{R}^{n\times 1} | e_1^T \cdot x \geq 0\} = \mathbb{R}^{n\times 1} \setminus \text{int}(H)$. Then the following inclusions are true
\begin{enumerate}
\item
\begin{align}\label{EA.3}
\mathcal{H} \cap \mathcal{D}_2 \subseteq \mathcal{H} \cap \mathcal{D}_1 
\end{align}
\item 
\begin{align}\label{EA.4}
\mathcal{G} \cap \mathcal{D}_1 \subseteq \mathcal{G} \cap \mathcal{D}_2
\end{align}
\end{enumerate}
\end{lemma}
\begin{proof}
Let $x = q \cdot e_1 + v$ with $e_1^T \cdot v = 0$ and assume that $x \in \mathcal{H} \cap \mathcal{D}_2$, i.e $q \leq 0$ and $\| q \cdot e_1 + v - (-q_2 \cdot e_1)\|^2 \leq r_2^2$. It follows that $(q+q_2)^2 + \|v\|^2 \leq r_2^2$. We want to check if $x \in \mathcal{D}_1$ i.e $\|q \cdot e_1 + v + q_1\cdot e_1\|^2 \leq^{?} r_1^2$ i.e $(q + q_1)^2 + \|v\|^2 \leq^{?} r_1^2$. However since
\begin{align}
(q + q_1)^2 + \|v\|^2 \leq (q + q_1)^2 + r_2^2 - (q + q_2)^2  
\end{align} we only prove that 
\begin{align}
(q + q_1)^2 + r_2^2 - (q + q_2)^2  \leq r_1^2 \iff (q + q_1)^2 - (q + q_2)^2 \leq r_1^2 - r_2^2
\end{align} that is 
\begin{align}\label{EA.6}
(q_1 - q_2)\cdot (2 \cdot q + (q_1 + q_2)) = 2 \cdot (q_1-q_2) \cdot q + q_1^2 - q_2^2 \leq r_1^2 - r_2^2
\end{align} But from (\ref{EA.2}) one has $_1^2 - r_2^2 = q_1^2 - q_2^2$ hence (\ref{EA.6}) is equivalent to
\begin{align}\label{EA.7}
2 \cdot (q_1 - q_2) \cdot q \leq 0
\end{align} But since $ |q_1| > |q_2|$ and $q_1 > 0$ follows that $q_1 - q_2 > 0$. Finally because $q \leq 0$ (\ref{EA.7}) is true and so is the claim (\ref{EA.3}).

 The claim in (\ref{EA.4}) is easily proved in a similar fashion. 
\end{proof}

Next we have 
\begin{lemma} Let us consider the following n-disks 
\begin{align}
\mathcal{D}_1 &= \{x \in \mathbb{R}^{n\times 1} | \|-q_1\cdot e_1 - x\| \leq r_1\} \nonumber \\
\mathcal{D}_2 &= \{x \in \mathbb{R}^{n\times 1}| \|-q_2 \cdot e_1 - x\| \leq r_2\} \nonumber \\
\mathcal{D}_3 &= \{x \in \mathbb{R}^{n\times 1} | \|-q_3 \cdot e_1 - x\| \leq r_3\}
\end{align} with $q_1, q_2 > 0$ and $r_1 > r_2 > r_3 \geq 0$ such that exists 
\begin{align}
0<a^2 = r_1^2 - q_1^2  = r_2^2 - q_2^2 = r_3^2 - q_3^2
\end{align} that is the n-disks share a common n-sphere i.e $\{x | e_1^T\cdot x = 0, \|x\| = a\}$. Let $\mathcal{H} = \{x \in \mathbb{R}^{n\times 1} | e_1^T \cdot x \leq 0\}$ and $\mathcal{G} = \{x \in \mathbb{R}^{n\times 1} | e_1^T \cdot x \geq 0\} = \mathbb{R}^{n\times 1} \setminus \text{int}(H)$. Then the following inclusions are true
\begin{enumerate}
\item
\begin{align}\label{EA.11}
\mathcal{H} \cap \mathcal{D}_3 \subseteq \mathcal{D}_1 \cap \mathcal{D}_3 \subseteq \mathcal{D}_2 \cap \mathcal{D}_3 \subseteq \mathcal{D}_2  
\end{align}
\item 
\begin{align}\label{EA.12}
\mathcal{G} \cap \mathcal{D}_1 \subseteq \mathcal{G} \cap \mathcal{D}_2 \subseteq \mathcal{G} \cap \mathcal{D}_3
\end{align}
\end{enumerate}
\end{lemma}
\begin{proof}
For (\ref{EA.12}) once can successively apply Lemma \ref{LA1} claim (\ref{EA.4}) to obtain the desired result. 

For (\ref{EA.11}) the last inclusion is obvious. For the first inclusion, let $x \in \mathcal{H}\cap \mathcal{D}_3$. Then we already have from Lemma \ref{LA1} that $x \in \mathcal{H} \cap \mathcal{D}_1 \cap \mathcal{D}_3 \subseteq \mathcal{D}_1 \cap \mathcal{D}_3$. For the second inclusion we prove the following:
\begin{align}
\mathcal{D}_1 \cap \mathcal{D}_3 \cap \mathcal{H} \subseteq \mathcal{D}_2 \cap \mathcal{D}_3 \cap \mathcal{H} \hspace{0.2cm}\text{and} \hspace{0.2cm}
\mathcal{D}_1 \cap \mathcal{D}_3 \cap \mathcal{G} \subseteq \mathcal{D}_2 \cap \mathcal{D}_3 \cap \mathcal{G}
\end{align} Indeed, in the above equation for the first inclusion let $x \in \mathcal{D}_1 \cap \mathcal{D}_3 \cap \mathcal{H}$. We have from Lemma \ref{LA1} that $\mathcal{D}_3 \cap \mathcal{H} \subseteq \mathcal{D}_2 \cap \mathcal{H}$ hence (since $x \in \mathcal{D}_3$) 
\begin{align}
x \in \mathcal{D}_1 \cap \mathcal{D}_2 \cap \mathcal{H} \cap \mathcal{D}_3 \subseteq \mathcal{D}_2  \cap \mathcal{D}_3 \cap \mathcal{H}
\end{align} Finally, let $x \in \mathcal{D}_1 \cap \mathcal{D}_3 \cap \mathcal{G}$. Applying Lemma \ref{LA1} we obtain that $ \mathcal{D}_1 \cap \mathcal{G} \subseteq  \mathcal{D}_2 \cap \mathcal{G}$. It follows 
\begin{align}
x \in  \mathcal{D}_3 \cap \mathcal{D}_2 \cap \mathcal{G} 
\end{align}
\end{proof}

\section{Approximating polytopes with a finite intersection of n-disks}
%
%

\begin{lemma}
Given $\delta > 0$ fixed, if $\exists \mathcal{B}(x,\delta) \subseteq \mathcal{P} $ then exists $\rho_{\delta} > 0$ such that for all $\rho \geq \rho_{\delta}$ one has 
\begin{align}
\mathcal{P} \subseteq \mathcal{C}_{\rho} \subseteq \bigcup_{x \in \mathcal{P}} \mathcal{B}(x,\delta)
\end{align}
\end{lemma}
\begin{proof} The first inclusion is true for $\rho \geq \bar{\rho}$ from Remark \ref{R3.5}. For the second inclusion one has the following. From Figure \ref{F1} one has
\begin{align}
x = r_1 - q_1 = \sqrt{q_1^2 + a^2} - q_1 = \frac{a^2}{q_1 + \sqrt{q_1^2 + a^2}}
\end{align} For a fixed $a$, letting $q_1$ be large enough one can see that $x \leq \delta$ will eventually occur. Therefore \textit{for any point $u$ in the intersection of the large disk with the right half-space exists a point $v$ in the intersection of the hyper-plane with the small disk such that $u \in \mathcal{B}(v,\delta)$}. 

\begin{figure}[h]
\centering
\includegraphics[scale=0.6]{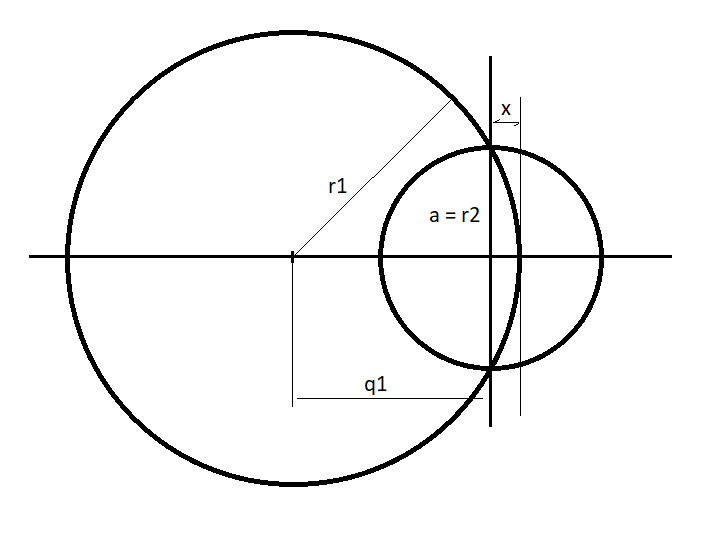}
\caption{Plane problem}
\label{F1}
\end{figure}
Consider w.l.o.g the set $\mathcal{H}_s = \{x | S^T \cdot x \leq 0\}$ i.e a facet of $\mathcal{P}$ . Since $\mathcal{C}_{\rho} \subseteq \bar{\mathcal{B}}\left( \frac{1}{2} \cdot 1_{n \times 1}, \frac{\sqrt{n}}{2} \right)$ we want to prove that 
\begin{align}
\bar{\mathcal{B}}(C_s, r_s) \cap \bar{\mathcal{B}}\left( \frac{1}{2} \cdot 1_{n \times 1}, \frac{\sqrt{n}}{2} \right)  \subseteq \bigcup_{x \in \mathcal{H}_s \cap \bar{\mathcal{B}}\left( \frac{1}{2} \cdot 1_{n \times 1}, \frac{\sqrt{n}}{2} \right) } \mathcal{B}(x, \delta)
\end{align} Regarding the above equation, it is easy to see that 
\begin{align}
\mathcal{H}_s \cap \bar{\mathcal{B}}(C_s, r_s) \cap \bar{\mathcal{B}}\left( \frac{1}{2} \cdot 1_{n \times 1}, \frac{\sqrt{n}}{2} \right)  & \subseteq \mathcal{H}_s \cap \bar{\mathcal{B}}\left( \frac{1}{2} \cdot 1_{n \times 1}, \frac{\sqrt{n}}{2} \right) \nonumber \\
& \subseteq \bigcup_{x \in \mathcal{H}_s \cap \bar{\mathcal{B}}\left( \frac{1}{2} \cdot 1_{n \times 1}, \frac{\sqrt{n}}{2} \right) } \mathcal{B}(x, \delta)
\end{align} Therefore we will now focus only on the elements of 
\begin{align}
\mathcal{G}_s \cap \bar{\mathcal{B}}(C_s, r_s) \cap \bar{\mathcal{B}}\left( \frac{1}{2} \cdot 1_{n \times 1}, \frac{\sqrt{n}}{2} \right) 
\end{align} where $\mathcal{G}_s = \{x | S^T \cdot x \geq 0\}$. For this in Figure \ref{F1}, take the large disk $\bar{\mathcal{B}}(C_s, r_s)$ and the small disk $\bar{\mathcal{B}}(P_s, \tilde{r}_s)$. Here consider $\|C_s - P_s\| = q_1$ and $a = \tilde{r_s}$. Then, as shown, given $\delta > 0$ exists $q_{\delta}$ such that for any $q_1 \geq q_{\delta}$ one has 
\begin{align} \label{EB.6}
\bar{\mathcal{B}}(C_s,r_s) \cap \mathcal{G}_s \subseteq \bigcup_{x \in \mathcal{H}_s \cap \mathcal{G}_s \cap \bar{\mathcal{B}}(P_s, \tilde{r}_s)} \mathcal{B}(x, \delta) 
\end{align} However, please note that due to construction one has
\begin{align}\label{EB.7}
\mathcal{H}_s \cap \mathcal{G}_s \cap \bar{\mathcal{B}}(P_s, \tilde{r}_s) &= \mathcal{H}_s \cap \mathcal{G}_s \cap \bar{\mathcal{B}}\left( \frac{1}{2} \cdot 1_{n \times 1}, \frac{1\sqrt{n}}{2} \right) \nonumber \\
& \subseteq \mathcal{H}_s \cap \bar{\mathcal{B}}\left( \frac{1}{2} \cdot 1_{n \times 1}, \frac{1\sqrt{n}}{2} \right)
\end{align} hence from (\ref{EB.6}) and (\ref{EB.7}) it is obtained:
\begin{align}
& \mathcal{G}_s \cap \bar{\mathcal{B}}(C_s, r_s) \cap \bar{\mathcal{B}}\left( \frac{1}{2} \cdot 1_{n \times 1}, \frac{\sqrt{n}}{2} \right) \subseteq \mathcal{G}_s \cap \bar{\mathcal{B}}(C_s, r_s)  \subseteq \nonumber \\
& \subseteq \bigcup_{x \in \mathcal{H}_s \cap \mathcal{G}_s \cap \bar{\mathcal{B}}\left( \frac{1}{2} \cdot 1_{n \times 1}, \frac{1\sqrt{n}}{2} \right)} \mathcal{B}(x, \delta) \nonumber \\
& \subseteq \bigcup_{x \in \mathcal{H}_s \cap \bar{\mathcal{B}}\left( \frac{1}{2} \cdot 1_{n \times 1}, \frac{1\sqrt{n}}{2} \right)} \mathcal{B}(x, \delta)
\end{align} 

Finally, since $\left\|P_s - \frac{1}{2} \cdot 1_{n\times 1} \right\|$ is fixed, form the existence of $q_{\delta}$ the existence of $\rho_{\delta}$ easily follows. 
\end{proof}

\section{Equivalent polytopes}
\begin{lemma} \label{LB.1}
For $R \geq R^{\star}$ exists $\hat{R}$ such that 
\begin{align}
\mathcal{P}_{R^2} =\hat{\mathcal{P}}_{\hat{R}^2}
\end{align} where
\begin{align}
\mathcal{P}_{R^2} &= \{x | h(x) - \|x - C\|^2 \leq -R^2\}  \nonumber \\ \hat{\mathcal{P}}_{\hat{R}^2} &= \{x | \hat{h}(x) - \|x - \hat{C}\|^2 \leq - \hat{R}^2\}
\end{align}
\end{lemma}
\begin{proof}Let $ R \geq R^{\star}$ and consider the inequalities:
\begin{align}\label{E3.61}
\|x - C_{k+}\|^2 - r_{k+}^2 - \|x - C\|^2 + R^2 \leq 0 \nonumber \\
\|x - \hat{C}_{k+}\|^2 - \hat{r}_{k+}^2 - \| x - \hat{C}\|^2 + \hat{R}^2 \leq 0
\end{align} We wonder if exists $\hat{R}$ such that the above inequalities are actually only one for all $k$. We have for the first inequality
\begin{align}
\left\|x - \left( \frac{1}{2}\cdot 1_{n \times 1} + q_k \cdot e_k \right)\right\|^2 - \left\| x - \left( \frac{1}{2} \cdot 1_{n \times 1} - \frac{\beta}{2} \cdot \frac{S}{\|S\|} \right) \right\|^2 - r_{k+}^2 + R^2 \leq 0
\end{align} which becomes the linear inequality
\begin{align}
&\left(\left( x - \left( \frac{1}{2}\cdot 1_{n \times 1} + q_k \cdot e_k \right)\right) - \left(x -  \left( \frac{1}{2} \cdot 1_{n \times 1} - \frac{\beta}{2} \cdot \frac{S}{\|S\|} \right) \right) \right)^T\cdot  \nonumber \\
& \cdot  \left( \left( x - \left( \frac{1}{2}\cdot 1_{n \times 1} + q_k \cdot e_k \right)\right) + \left(x -  \left( \frac{1}{2} \cdot 1_{n \times 1} - \frac{\beta}{2} \cdot \frac{S}{\|S\|} \right) \right) \right) - r_{k+}^2 + R^2 \leq 0
\end{align} then 
\begin{align}\label{E3.64}
\left( -q_k \cdot e_k - \frac{\beta}{2} \cdot \frac{S}{\|S\|}\right)^T \cdot \left( 2 \cdot x - 1_{n\times 1} - q_k\cdot e_k + \frac{\beta}{2}\cdot \frac{S}{\|S\|} \right) - r_{k+}^2 + R^2 \leq 0
\end{align} After similar calculations, the second inequality in (\ref{E3.61}) becomes
\begin{align} \label{E3.65}
\left( -\alpha \cdot q_k \cdot e_k - \frac{\alpha \cdot \beta}{2} \cdot \frac{S}{\|S\|}\right)^T \cdot \left( 2 \cdot x -1_{n \times 1}- \alpha \cdot q_k \cdot e_k + \frac{\alpha \cdot \beta}{2}\cdot \frac{S}{\|S\|} \right) - \hat{r}_{k+}^2 + \hat{R}^2 \leq 0
\end{align} 
 We would want to prove that exists $\hat{R}$ not depending on $k$ such that  the inequality in (\ref{E3.64}) multiplied by $\alpha$ becomes the inequality from (\ref{E3.65}). 

About $r_{k+}$ and $\hat{r}_{k+}$ we know from (\ref{E3.20}) and (\ref{E3.53}) that
\begin{align} \label{E3.66}
&\left\| \frac{1}{2} \cdot 1_{n \times 1} - \frac{1}{2} \cdot e_k - \left( \frac{1}{2}\cdot 1_{n\times 1} + q_k \cdot e_k\right)\right\|^2 + \left( \frac{\sqrt{n-1}}{2}\right)^2 = r_{k+1}^2 \nonumber \\
&\left\| \frac{1}{2} \cdot 1_{n \times 1} - \frac{1}{2} \cdot e_k - \left( \frac{1}{2}\cdot 1_{n\times 1} + \alpha \cdot q_k \cdot e_k\right)\right\|^2 + \left( \frac{\sqrt{n-1}}{2}\right)^2 = \hat{r}_{k+1}^2
\end{align} The first equality in (\ref{E3.66}) becomes:
\begin{align}\label{E3.67}
\left( \frac{1}{2} + q_k\right)^2 + \left( \frac{\sqrt{n-1}}{2}\right)^2 = r_{k+}^2
\end{align} while the second one becomes:
\begin{align} \label{E3.68}
\hat{r}_{k+}^2 &= \left( \frac{1}{2} + \alpha \cdot q_k\right)^2 + \left( \frac{\sqrt{n-1}}{2}\right)^2 
\end{align}

We write (\ref{E3.65}) as follows
\begin{align}
&\alpha \cdot \left( -q_k \cdot e_k - \frac{\beta}{2} \cdot \frac{S}{\|S\|}\right)^T \cdot \left( 2 \cdot x - 1_{n\times 1}- q_k \cdot e_k + \frac{\beta}{2}\cdot \frac{S}{\|S\|} \right) - \alpha \cdot \left(r_{k+}^2 - R^2\right) +  \nonumber \\
& \alpha \cdot \left( -q_k \cdot e_k - \frac{\beta}{2} \cdot \frac{S}{\|S\|}\right)^T \cdot (\alpha-1) \cdot \left( - q_k \cdot e_k + \frac{\beta}{2}\cdot \frac{S}{\|S\|} \right)  + \alpha \cdot \left(r_{k+}^2 - R^2\right) + \nonumber \\
& - \hat{r}_{k+}^2 + \hat{R}^2 \leq 0
\end{align} Let 
\begin{align}\label{E3.70}
&\theta_{k+} = \nonumber \\
& \alpha \left( -q_k \cdot e_k - \frac{\beta}{2} \cdot \frac{S}{\|S\|}\right)^T (\alpha-1) \left( - q_k \cdot e_k + \frac{\beta}{2}\cdot \frac{S}{\|S\|} \right)  + \alpha \left(r_{k+}^2 - R^2\right) - \hat{r}_{k+}^2 + \hat{R}^2
\end{align} We want to prove that exists $\hat{R}$ not depending on $k$ such that $\theta_{k+} = 0$. Indeed, we rewrite $\theta_{k+}$ as follows
\begin{align}\label{E3.71}
\alpha \cdot (\alpha - 1) \cdot \left( q_k^2 - \frac{\beta^2}{4} \right) + \alpha\cdot r_{k+}^2 - \hat{r}_{k+}^2 - \alpha \cdot R^2 + \hat{R}^2
\end{align} Finally we show that $\alpha\cdot (\alpha - 1) \cdot q_k^2 + \alpha\cdot r_{k+}^2 - \hat{r}_{k+}^2$ does not depend on $k$. Using (\ref{E3.67}) and  (\ref{E3.68}) we obtain
\begin{align}\label{E3.74}
\alpha\cdot r_{k+}^2-\hat{r}_{k+}^2& = \alpha \cdot \left( \left( \frac{1}{2} + q_k\right)^2 + \frac{n-1}{4} \right) - \left( \left(\frac{1}{2} + \alpha \cdot q_k \right)^2 + \frac{n-1}{4}\right) \nonumber \\
& = \alpha \cdot \left( \frac{n}{4} + q_k + q_k^2\right) - \left( \frac{n}{4} + \alpha\cdot q_k + \alpha^2 \cdot q_k^2 \right) \nonumber \\
& = (\alpha - 1) \cdot \frac{n}{4} + \alpha \cdot (1 - \alpha) \cdot q_k^2
\end{align}  hence 
\begin{align} \label{E3.76}
\theta_{k+} = (\alpha - 1) \cdot \frac{n}{4} + \alpha \cdot (\alpha - 1) \frac{-\beta^2}{4} - \alpha \cdot R^2 + \hat{R}^2
\end{align}
We will not make the calculations for what would be $\theta_{k-}$ and simply assume (because of the high similarity) that $\theta_{k-} = \theta_{k+}$. 

Let us now consider the inequalities
\begin{align}
&\|x - C_h\|^2 - r_h^2 - \| x - C\|^2 + R^2 \leq 0 \nonumber \\
& \|x - \hat{C}_h\|^2 - \hat{r}_h^2 - \| x - \hat{C}\|^2 + \hat{R}^2 \leq 0
\end{align} The first of them becomes:
\begin{align}
\left\| x - \left( \frac{1}{2} \cdot 1_{n \times 1} + q_h \cdot \frac{1_{n\times 1}}{\|1_{n \times 1}\|}\right) \right\|^2 - \left\| x - \left(\frac{1}{2} \cdot 1_{n \times 1} - \frac{\beta}{2} \cdot \frac{S}{\|S\|} \right) \right\|^2 - r_h^2 + R^2 \leq 0
\end{align} then 
\begin{align}
\left( \frac{-\beta}{2} \frac{S}{\|S\|} - q_h \cdot \frac{1_{n \times 1}}{1_{n \times 1}} \right)^T \cdot \left( 2 \cdot x -1_{n \times 1} - q_h\cdot \frac{1_{n \times 1}}{\|1_{n \times 1}\|} + \frac{\beta}{2} \cdot \frac{S}{\|S\|} \right) -r_h^2 + R^2 \leq 0
\end{align} Similarly, the second inequality becomes:
\begin{align}\label{E3.79}
\alpha \cdot \left( \frac{-\beta}{2} \frac{S}{\|S\|} - q_h  \frac{1_{n \times 1}}{1_{n \times 1}} \right)^T \cdot \left( 2 \cdot x -1_{n \times 1} -  \alpha q_h \frac{1_{n \times 1}}{\|1_{n \times 1}\|} + \frac{\alpha  \beta}{2}  \frac{S}{\|S\|} \right) -\hat{r}_h^2 + \hat{R}^2 \leq 0
\end{align} From (\ref{E3.29}) and (\ref{E3.53}) we obtain:
\begin{align}
&\| C_h - P_h\|^2 + \tilde{r}_h^2 = r_h^2 = \left\| \frac{1}{2} \cdot 1_{n \times 1} + q_h \cdot \frac{1_{n \times 1}}{\|1_{n \times 1}\|} - P_h \right\|^2 + \tilde{r}_h^2 \nonumber \\
&\| \hat{C}_h - P_h\|^2 + \tilde{r}_h^2 = \hat{r}_h^2 = \left\|\frac{1}{2} \cdot 1_{n \times 1} + \alpha \cdot q_h \cdot \frac{1_{n \times 1}}{\|1_{n \times 1}\|} - P_h \right\|^2 + \tilde{r}_h^2
\end{align}  We rewrite (\ref{E3.79}) as follows
\begin{align}
&\alpha \cdot \left( \frac{-\beta}{2} \frac{S}{\|S\|} - q_h  \frac{1_{n \times 1}}{1_{n \times 1}} \right)^T \cdot \left( 2 \cdot x -1_{n \times 1} -  q_h \frac{1_{n \times 1}}{\|1_{n \times 1}\|} + \frac{\beta}{2}  \frac{S}{\|S\|} \right) + \alpha (-r_h^2 + R^2) + \nonumber \\
&\alpha \cdot \left( \frac{-\beta}{2} \frac{S}{\|S\|} - q_h  \frac{1_{n \times 1}}{1_{n \times 1}} \right)^T \cdot (\alpha - 1) \left( -q_h\cdot \frac{1_{n\times 1}}{\|1_{n \times 1}\|} + \frac{\beta}{2} \frac{S}{\|S\|} \right) - \alpha (-r_h^2 + R^2) - \nonumber \\
& - \hat{r}_h^2 + \hat{R}^2 \leq 0
\end{align} Let us denote
\begin{align}
\theta_h  = \alpha \cdot (\alpha - 1) \cdot \left( q_h^2 - \frac{\beta^2}{4}\right) + \alpha \cdot r_h^2 - \hat{r}_h^2 - \alpha \cdot R^2 + \hat{R}^2
\end{align} Here we evaluate:
\begin{align}
\alpha \cdot r_h^2 - &\hat{r}_h^2 =  \alpha \cdot \left( \left\|\frac{1}{2} \cdot 1_{n \times 1} + q_h \cdot \frac{1_{n \times 1}}{\|1_{n \times 1}\|} - P_h \right\|^2 + \tilde{r}_h^2\right) - \nonumber \\ & \left( \left\|\frac{1}{2} \cdot 1_{n \times 1} + \alpha \cdot q_h \cdot \frac{1_{n \times 1}}{\|1_{n \times 1}\|} - P_h \right\|^2 + \tilde{r}_h^2 \right) \nonumber \\
= \alpha \cdot &\left( \left\| \left(\frac{1}{2} \cdot 1_{n \times 1} - P_h\right) + q_h \cdot \frac{1_{n\times 1}}{\|1_{n \times 1}\|}\right\|^2 + \tilde{r}_h^2 \right) - 
\nonumber \\
& \left( \left\| \left(\frac{1}{2} \cdot 1_{n \times 1} - P_h\right) + \alpha\cdot q_h \cdot \frac{1_{n\times 1}}{\|1_{n \times 1}\|}\right\|^2 + \tilde{r}_h^2 \right) \nonumber \\
= \alpha \cdot &\left(\left\|\frac{1}{2} \cdot 1_{n \times 1} - P_h \right\|^2 + \tilde{r}_h^2 + q_h^2 + 2 \cdot q_h \cdot \left( \frac{1}{2} \cdot 1_{n \times 1}-P_h \right) ^T \cdot \frac{1_{n\times 1}}{\|1_{n \times 1}\|}\right) - \nonumber \\
& \left(\left\|\frac{1}{2} \cdot 1_{n \times 1} - P_h \right\|^2 + \tilde{r}_h^2 + \alpha^2 \cdot q_h^2 + 2 \cdot \alpha \cdot q_h \cdot \left( \frac{1}{2} \cdot 1_{n \times 1}-P_h \right) ^T \cdot \frac{1_{n\times 1}}{\|1_{n \times 1}\|}\right)
\end{align} But from (\ref{E3.28}) one has $\frac{n}{4} = \left\| \frac{1}{2} \cdot 1_{n \times 1} - P_h\right\|^2 + \tilde{r}_h^2$ hence 
\begin{align}
\alpha\cdot r_h^2 - \hat{r}_h^2 = (\alpha - 1) \cdot \frac{n}{4} + \alpha \cdot q_h^2 - \alpha^2 \cdot q_h^2 
\end{align} and 
\begin{align}\label{E3.85}
\theta_h = \alpha \cdot (\alpha - 1) \cdot \frac{-\beta^2}{4} + (\alpha-1) \cdot \frac{n}{4} - \alpha \cdot R^2 + \hat{R}^2
\end{align}

Finally, we consider the inequalities:
\begin{align} 
&\|x - C_s\|^2 - r_s^2 - \|x - C\|^2 + R^2 \leq 0\nonumber \\
&\|x - \hat{C}_s\|^2 - \hat{r}_s^2 - \|x - \hat{C}\|^2 + \hat{R}^2 \leq 0
\end{align} The first inequality becomes:
\begin{align}
\left\|x - \left( \frac{1}{2} \cdot 1_{n \times 1} -q_s \cdot \frac{S}{\|S\|} \right)\right\|^2 - \left\|x - \left( \frac{1}{2} \cdot 1_{n \times 1} -\frac{\beta}{2} \cdot \frac{S}{\|S\|} \right)\right\|^2 - r_s^2 + R^2 \leq 0
\end{align} then 

\begin{align}
\left( q_s\frac{S}{\|S\|} - \frac{\beta}{2} \cdot \frac{S}{\|S\|}\right)^T \cdot \left( 2 \cdot x - 1_{n\times 1} + q_s\frac{S}{\|S\|} + \frac{\beta}{2} \frac{S}{\|S\|} \right) - r_s^2 + R^2 \leq 0
\end{align} and similarly the second inequality becomes:
\begin{align}\label{E3.88}
\alpha \cdot \left( q_s\frac{S}{\|S\|} - \frac{\beta}{2} \cdot \frac{S}{\|S\|}\right)^T \cdot \left( 2 \cdot x - 1_{n\times 1} + \alpha q_s\frac{S}{\|S\|} + \frac{\alpha\beta}{2} \frac{S}{\|S\|} \right) - \hat{r}_s^2 + \hat{R}^2 \leq 0
\end{align} 
From (\ref{E3.25}) and (\ref{E3.53}) we obtain:
\begin{align}
&\| C_s - P_s\|^2 + \tilde{r}_s^2 = r_s^2 = \left\| \frac{1}{2} \cdot 1_{n \times 1} - q_s \cdot \frac{S}{\|S\|} - P_s \right\|^2 + \tilde{r}_s^2 \nonumber \\
&\| \hat{C}_s - P_s\|^2 + \tilde{r}_s^2 = \hat{r}_s^2 = \left\|\frac{1}{2} \cdot 1_{n \times 1} - \alpha \cdot q_s \cdot \frac{S}{\|S\|} - P_s \right\|^2 + \tilde{r}_s^2
\end{align}

As we did above with the other so we do here. We focus on (\ref{E3.88}) and rewrite it as
\begin{align}
&\alpha \cdot \left( q_s\frac{S}{\|S\|} - \frac{\beta}{2} \cdot \frac{S}{\|S\|}\right)^T \cdot \left( 2 \cdot x - 1_{n\times 1} + q_s\frac{S}{\|S\|} + \frac{\beta}{2} \frac{S}{\|S\|} \right) + \alpha \cdot \left( -r_s^2 + R^2 \right) + \nonumber \\
&\alpha \cdot \left( q_s\frac{S}{\|S\|} - \frac{\beta}{2} \cdot \frac{S}{\|S\|}\right)^T \cdot (\alpha-1) \cdot \left( q_s\frac{S}{\|S\|} + \frac{\beta}{2} \frac{S}{\|S\|} \right) - \alpha \cdot (-r_s^2 + R^2) -\nonumber \\
&-\hat{r}_s^2 + \hat{R}^2 \leq 0
\end{align} Let 
\begin{align}
\theta_s = \alpha \cdot (\alpha - 1) \cdot \left( q_s^2 - \frac{\beta^2}{4}\right) + \alpha \cdot r_s^2 - \hat{r}_s^2 - \alpha \cdot R^2 + \hat{R}^2
\end{align} and we evaluate
\begin{align}
\alpha \cdot r_s^2 - &\hat{r}_s^2 =  \alpha \cdot \left( \left\|\frac{1}{2} \cdot 1_{n \times 1} - q_s \cdot \frac{S}{\|S\|} - P_s \right\|^2 + \tilde{r}_s^2\right) - \nonumber \\ & \left( \left\|\frac{1}{2} \cdot 1_{n \times 1} - \alpha \cdot q_s \cdot \frac{S}{\|S\|} - P_s \right\|^2 + \tilde{r}_s^2 \right) \nonumber \\
= \alpha \cdot &\left( \left\| \left(\frac{1}{2} \cdot 1_{n \times 1} - P_s\right) - q_s \cdot \frac{S}{\|S\|}\right\|^2 + \tilde{r}_s^2 \right) - 
\nonumber \\
& \left( \left\| \left(\frac{1}{2} \cdot 1_{n \times 1} - P_s\right) - \alpha\cdot q_s \cdot \frac{S}{\|S\|}\right\|^2 + \tilde{r}_s^2 \right) \nonumber \\
= \alpha \cdot &\left(\left\|\frac{1}{2} \cdot 1_{n \times 1} - P_s \right\|^2 + \tilde{r}_s^2 + q_s^2 - 2 \cdot q_s \cdot \left( \frac{1}{2} \cdot 1_{n \times 1}-P_s \right) ^T \cdot \frac{S}{\|S\|}\right) - \nonumber \\
& \left(\left\|\frac{1}{2} \cdot 1_{n \times 1} - P_s \right\|^2 + \tilde{r}_s^2 + \alpha^2 \cdot q_s^2 - 2 \cdot \alpha \cdot q_s \cdot \left( \frac{1}{2} \cdot 1_{n \times 1}-P_s \right) ^T \cdot \frac{S}{\|S\|}\right)
\end{align}  But from (\ref{E3.24}) one has $\frac{n}{4} = \left\| \frac{1}{2} \cdot 1_{n \times 1} - P_s\right\|^2 + \tilde{r}_s^2$ hence 
\begin{align}
\alpha\cdot r_s^2 - \hat{r}_s^2 = (\alpha - 1) \cdot \frac{n}{4} + \alpha \cdot q_s^2 - \alpha^2 \cdot q_s^2 
\end{align} and 
\begin{align}\label{E3.93}
\theta_s = \alpha \cdot (\alpha - 1) \cdot \frac{-\beta^2}{4} + (\alpha-1) \cdot \frac{n}{4} - \alpha \cdot R^2 + \hat{R}^2
\end{align}

For $R = R^{\star}$, from (\ref{E3.76}), (\ref{E3.85}) and (\ref{E3.93}) in order to have $\theta_{k+} = \theta_{k-} = \theta_h = \theta_s = 0$ we take 
\begin{align}
\hat{R}^2 = \alpha \cdot R^2 + \alpha \cdot (\alpha-1) \cdot \frac{\beta^2}{4} - (\alpha - 1)\cdot \frac{n}{4} 
\end{align} 

\end{proof}

\end{document}